\documentclass[letterpaper, 10pt, conference]{ieeeconf} \IEEEoverridecommandlockouts 
\usepackage{setspace, epsfig, psfrag, amssymb, amsmath, cite, color} \newfont{\smalll}{cmr8}

\def\IR{\mathbb{R}}
\def\IS{\mathbb{S}}
\def\zero{\mathbf{0}}
\def\one{\mathbf{1}}

\def\diag{\mathrm{diag}}
\def\span{\mathrm{span}}
\def\eye{\mathrm{I}}

\def\IS{\hbox{I\hskip-.1em S}}
\def\IC{\hbox{C\hskip-
.5em\raise.5ex\hbox{$\scriptscriptstyle\mid$}}\ }
\def\Ic{\hbox{\smalll C\hskip-
.5em\raise.3ex\hbox{$\scriptscriptstyle\mid$}}\ }

\def\T={\buildrel {\scriptscriptstyle\triangle} \over =}

\def\sqr#1#2{{\vcenter{\vbox{\hrule height.#2pt\hbox{\vrule
width.#2pt height#1pt \kern#1pt\vrule width.#2pt}\hrule
height.#2pt}}}}
\def\square{\mathchoice\sqr64\sqr64\sqr33\sqr33}
\def\diag{\mathop{\rm diag}}

\def\block-diag{\mathop{\rm block{\scriptstyle -}diag}}

\def\pmbb#1{\setbox0=\hbox{#1}\raise 0.5ex\box0}

\newcommand{\bequ}{\begin{eqnarray}}
\newcommand{\eequ}{\end{eqnarray}}
\newcommand{\mT}{^\mathrm{T}}

\newcommand{\rom}{\mathrm}

\newcommand {\teq}      {\triangleq}

\newcommand {\beq}      {\begin{equation}}
\newcommand {\eeq}      {\end{equation}}

\def\IR{{\mathbb R}}
\def\IC{{\mathbb C}}

\def\IS{{\mathbb S}}

\renewcommand\theenumi{\roman{enumi}}
\renewcommand\theenumii{\alph{enumii}}
\renewcommand\theenumiii{\arabic{enumiii}}
\renewcommand\theenumiv{\Alph{enumiv}}

\def\sc{s_{\rm c}}

\def\sc{s_{\rm c}}

\begin{document}
    {
    \title{{  
    {\bf On Networks with Active and Passive Agents}   
    }}
    }{
\author{Tansel Yucelen
\thanks{T. Yucelen is an Assistant Professor of the Mechanical and Aerospace Engineering Department and the Director of the Advanced Systems Research Laboratory at the Missouri University of Science and Technology, Rolla, MO 65409, USA (e-mail: {\tt\small yucelen@mst.edu}).
}} 
\markboth{} {Shell \MakeLowercase{\textit{et al.}}: Bare Demo of IEEEtran.cls for Journals} \newcommand{\eqnref}[1]{(\ref{#1})}
\newcommand{\class}[1]{\texttt{#1}} \newcommand{\package}[1]{\texttt{#1}} \newcommand{\file}[1]{\texttt{#1}} \newcommand{\BibTeX}{\textsc{Bib}\TeX}
\maketitle 

\onehalfspacing


\begin{abstract}

We introduce an active--passive networked multiagent system framework, which consists of agents subject to exogenous inputs (active agents) and agents without any inputs (passive agents), and analyze its convergence using Lyapunov stability.

\end{abstract}


\subsection{Preliminaries}

In the multiagent literature, graphs are broadly adopted to encode interactions in networked systems \cite{ref:1,ref:2}.
An \textit{undirected} graph $\mathcal{G}$ is defined by a set $\mathcal{V}_\mathcal{G}=\{1,\ldots,n\}$ of \textit{nodes}
and a set $\mathcal{E}_\mathcal{G} \subset \mathcal{V}_\mathcal{G} \times \mathcal{V}_\mathcal{G}$ of \textit{edges}.
If $(i,j) \in \mathcal{E}_\mathcal{G}$, then the nodes $i$ and $j$ are \textit{neighbors} and the neighboring relation is indicated with $i \sim j$.
The \textit{degree} of a node is given by the number of its neighbors.
Letting $d_i$ be the degree of node $i$, then the \textit{degree} matrix of a graph $\mathcal{G}$, $\mathcal{D}(\mathcal{G}) \in \IR^{n \times n}$, is given by $\mathcal{D}(\mathcal{G}) \triangleq \diag(d)$, $d=[d_1,\ldots,d_n]\mT$. 
A \textit{path} $i_0 i_1 \ldots i_L$ is a finite sequence of nodes such that $i_{k-1} \sim i_k$, $k=1, \ldots, L$,
and a graph $\mathcal{G}$ is \textit{connected} if there is a path between any pair of distinct nodes.
The \textit{adjacency} matrix of a graph $\mathcal{G}$, $\mathcal{A}(\mathcal{G}) \in \IR^{n \times n}$, is given by
\begin{eqnarray}
   [\mathcal{A}(\mathcal{G})]_{ij}
   &\teq&
   \left\{ \begin{array}{cl}
      1, & \mbox{ if $(i,j)\in\mathcal{E}_\mathcal{G}$},
      \\
      0, & \mbox{otherwise}.
   \end{array} \right.
   \label{AdjMat}
\end{eqnarray}
The \textit{Laplacian} matrix of a graph, $\mathcal{L}(\mathcal{G}) \in \overline{\IS}_+^{\hspace{0.1em} n \times n}$, playing a central role in many graph theoretic treatments of multiagent systems, is given by
\bequ
    \mathcal{L}(\mathcal{G}) &\triangleq& \mathcal{D}(\mathcal{G}) - \mathcal{A}(\mathcal{G}). \label{Laplacian}
\eequ
Throughout this note, we model a given multiagent system by a connected, undirected graph $\mathcal{G}$, where nodes and edges represent agents and inter-agent communication links, respectively. 


\vspace{0.25cm}

\subsection{Problem Formulation}

Consider a system of $n$ agents exchanging information among each other using their local measurements according to a connected, undirected graph $\mathcal{G}$. 
In addition, consider that there exists $m\ge1$ exogenous inputs that interact with this system. 
We make the following definitions.

\textbf{Definition 1.} If agent $i$, $i=1,\ldots,n$, is subject to one or more exogenous inputs (resp., no exogenous inputs), then it is an active agent (resp., passive agent).

\textbf{Definition 2.} If an exogenous input interacts with only one agent (resp., multiple agents), then it is an isolated input (resp., non-isolated input).

In this note, we are interested in the problem of driving the states of all (active and passive) agents to the average of the applied exogenous inputs. 
Motivating from this standpoint, we propose an integral action-based distributed control approach given by
\vspace{-0.05cm}
\bequ
	\dot{x}_i(t)&=&-\sum_{i \sim j}\Bigl(x_i(t)-x_j(t)\Bigl)+\sum_{i \sim j}\Bigl(\xi_i(t)-\xi_j(t)\Bigl) \nonumber\\
		        &&-\sum_{i \sim h}\Bigl(x_i(t)-c_h\Bigl), \quad x_i(0)=x_{i0}, \label{eq:1} \\
	\dot{\xi}_i(t) &=& -\sum_{i \sim j}\Bigl(x_i(t)-x_j(t)\Bigl), \quad \xi_i(0)=\xi_{i0}, \label{eq:2}	        
\eequ
where $x_i(t)\in\IR$ and $\xi_i(t)\in\IR$ denote the state and the integral action of agent $i$, $i=1,\ldots,n$, respectively, 
and $c_h\in\IR$, $h=1,\ldots,m$, denotes an exogenous input applied to this agent. 
Similar to the $i \sim j$ notation indicating the neighboring relation between agents, we use $i \sim h$ to indicate the exogenous inputs that an agent is subject to. 

Next, let $x(t)=\bigl[x_1(t),x_2(t),\ldots,x_n(t)\bigl]\mT\in\IR^n$, $\xi(t)=\bigl[\xi_1(t),\xi_2(t),\ldots,\xi_n(t)\bigl]\mT\in\IR^n$, and $c=\bigl[c_1,c_2,$ $\ldots,c_m,0,\ldots,0\bigl]\in\IR^n$, where we assume $m \le n$ for the ease of the following notation and without loss of generality.
We can now write (\ref{eq:1}) and (\ref{eq:2}) in a compact form as
\vspace{0.1cm}
\bequ
	\dot{x}(t)&=&-\mathcal{L}(\mathcal{G})x(t)+\mathcal{L}(\mathcal{G})\xi(t)-K_1x(t)+K_2c, \nonumber\\ 
	&& \hspace{4.2cm} x(0)=x_0, \label{eq:10} \\ 
	\dot{\xi}(t)&=&-\mathcal{L}(\mathcal{G})x(t), \quad \xi(0)=\xi_0, \label{eq:11}
\eequ
where $\mathcal{L}(\mathcal{G}) \in \overline{\IS}_+^{\hspace{0.1em} n \times n}$, 
\vspace{0.1cm}
\bequ
	K_1&\triangleq&\rom{diag}([k_{1,1},k_{1,2},\ldots,k_{1,n}]\mT)\in\overline{\IS}_+^{\hspace{0.1em} n \times n},
\eequ 
with $k_{1,i}\in\overline{\mathbb{Z}}_+$ denoting the number of the exogenous inputs applied to agent $i$, $i=1,\ldots,n$, and
\vspace{0.1cm}
\bequ
 K_2&\triangleq&\begin{bmatrix}
 k_{2,11} & k_{2,12}& \cdots & k_{2,1n}\\
 k_{2,21} & k_{2,22}& \cdots & k_{2,2n} \\
  \vdots & \vdots & \ddots &\vdots \\
  k_{2,n1}& k_{2,n2}&\cdots & k_{2,nn}
 \end{bmatrix}\in\IR^{n \times n}, 
\eequ
with $k_{2,ih}=1$ if the exogenous input $c_h(t)$, $h=1,\ldots,m$, is applied to agent $i$, $i=1,\ldots,n$, and $k_{2,ih}=0$ otherwise. 
Note that $k_{1,i}=\sum_{j=1}^{n}k_{2,ij}$.

Since we are interested in driving the states of all (active and passive) agents to the average of the applied exogenous inputs, 
let
\vspace{0.1cm}
\bequ
	\delta(t) &\triangleq& x(t) - \epsilon\one_n\in\IR^n, \label{delta:error} \\
	\epsilon &\triangleq& \frac{\one_n\mT K_2 c}{\one_n\mT K_2 \one_n}\in\IR, \label{epsilon:e}
\eequ
be the error between $x_i(t)$, $i=1,\ldots,n$, and the average of the applied exogenous inputs $\epsilon$. 
Based on (\ref{epsilon:e}), $\epsilon$ can be equivalently written as
\vspace{0cm}
\bequ
\epsilon&=&\Bigl(k_{2,11}c_1+k_{2,12}c_2+\cdots+k_{2,21}c_1\nonumber\\&&+k_{2,22}c_2+\cdots\Bigl)/\Bigl(k_{2,11}+k_{2,12}\nonumber\\&&+\cdots+k_{2,21}+k_{2,21}+\cdots\Bigl), \label{hey:1}
\eequ 
which is the average of the applied exogenous inputs. 

\vspace{0.25cm}

\subsection{Convergence Analysis}

It follows from (\ref{delta:error}) and $\mathcal{L}(\mathcal{G}) \one_n = \zero_n$ of Lemma 1 that
\bequ
	\dot{\delta}(t) &=& -\mathcal{L}(\mathcal{G})\bigl[\delta(t)+\epsilon\one_n\bigl]+\mathcal{L}(\mathcal{G})\xi(t)-K_1\bigl[\delta(t)\nonumber\\&&+\epsilon\one_n\bigl]+K_2c(t) \nonumber\\
		&=& -\mathcal{F}(\mathcal{G}) \delta(t)+\mathcal{L}(\mathcal{G})\xi(t)-\bigl[ K_1\one_n \epsilon -K_2 c\bigl] \nonumber\\
		&=& -\mathcal{F}(\mathcal{G}) \delta(t)+\mathcal{L}(\mathcal{G})\xi(t)-\biggl[\frac{K_1 \one_n \one_n\mT K_2 c}{\one_n\mT K_2 \one_n} - K_2 c \biggl] \nonumber\\
		&=& -\mathcal{F}(\mathcal{G}) \delta(t)+\mathcal{L}(\mathcal{G})\xi(t)-L_cK_2c, \label{it:1}
\eequ
where $\mathcal{F}(\mathcal{G})\triangleq\mathcal{L}(\mathcal{G})+K_1$ and 
\bequ
	L_c &\triangleq& \frac{K_1 \one_n \one_n\mT}{\one_n\mT K_2 \one_n} - \eye_n. \label{it:2}
\eequ
Note that $\mathcal{F}(\mathcal{G})\in\IS_+^{n \times n}$ and 
\bequ
	\one_n\mT L_c = \one_n\mT\bigg[\frac{K_1 \one_n \one_n\mT}{\one_n\mT K_2 \one_n} - \eye_n\biggl] 
	 = \frac{\one_n\mT K_1 \one_n }{\one_n\mT K_2 \one_n}\one_n\mT - \one_n\mT 
	 = 0, \label{extra:e1}
\eequ
since $(\one_n\mT K_1 \one_n )/(\one_n\mT K_2 \one_n)=1$ from $k_{1,i}=\sum_{j=1}^{n}k_{2,ij}$.

Next, letting
\bequ
	e(t) &\triangleq& \xi(t)-\mathcal{L}^\dagger (\mathcal{G}) L_c K_2 c, \label{it:3}
\eequ
and using (\ref{it:3}) in (\ref{it:1}) yields
\vspace{0cm}
\bequ
	\dot{\delta}(t) &\hspace{-0.15cm}=\hspace{-0.15cm}& -\mathcal{F}(\mathcal{G}) \delta(t)+\mathcal{L}(\mathcal{G})\bigl[e(t)+\mathcal{L}^\dagger (\mathcal{G}) L_c K_2 c\bigl]-L_cK_2c \nonumber\\
	&\hspace{-0.15cm}=\hspace{-0.15cm}&  -\mathcal{F}(\mathcal{G}) \delta(t)+\mathcal{L}(\mathcal{G})e(t)+\Bigl[\eye_n-\frac{1}{n}\one_n\one_n\mT\Bigl] L_c K_2 c\nonumber\\&& -L_cK_2c \nonumber\\
	&\hspace{-0.15cm}=\hspace{-0.15cm}&  -\mathcal{F}(\mathcal{G}) \delta(t)+\mathcal{L}(\mathcal{G})e(t), \label{it:444}
\eequ
since $\frac{1}{n}\one_n\one_n\mT L_c K_2 c=0$ as a direct consequence of (\ref{extra:e1}). 
In addition, differentiating (\ref{it:3}) with respect to time yields
\bequ
	\dot{e}(t)&=& -\mathcal{L}(\mathcal{G})\bigl[\delta(t)+\epsilon \one_n \bigl] \nonumber\\
	&=& -\mathcal{L}(\mathcal{G})\delta(t), \label{it:555}
\eequ
where $\mathcal{L}(\mathcal{G}) \one_n = \zero_n$. 
The following theorem shows that the state of all agents $x_i(t)$, $i=1,\ldots,n$ asymptotically converge to $\epsilon$.

\textbf{Theorem 1.} Consider the networked multiagent system given by (\ref{eq:1}) and (\ref{eq:2}), where agents exchange information using local measurements and with $\mathcal{G}$ defining a connected, undirected graph topology. 
Then, the closed-loop error dynamics defined by (\ref{it:444}) and (\ref{it:555}) are Lyapunov stable for all initial conditions and $\delta(t)$ asymptotically vanishes.

\textit{Proof.} Proof follows by considering Lyapunov function candidate given by $V(\delta,e)=\frac{1}{2}\delta\mT\delta+\frac{1}{2}e\mT e$ and differentiating it along the trajectories of (\ref{it:444}) and (\ref{it:555}). 
\hfill $\square$

Note that a generalized version of the proposed integral action-based distributed control approach can be given by
\bequ
	\dot{x}_i(t)&=&-\alpha\sum_{i \sim j}\Bigl(x_i(t)-x_j(t)\Bigl)+\sum_{i \sim j}\Bigl(\xi_i(t)-\xi_j(t)\Bigl) \nonumber\\
		        &&-\alpha \sum_{i \sim h}\Bigl(x_i(t)-c_h\Bigl), \quad x_i(0)=x_{i0}, \label{eq:199} \\
	\dot{\xi}_i(t) &=& -\gamma \sum_{i \sim j}\Bigl(x_i(t)-x_j(t)\Bigl), \quad \xi_i(0)=\xi_{i0}, \label{eq:299}	        
\eequ
where $\alpha \in \IR_+$ and $\gamma \in \IR_+$. 

\vspace{0.25cm}

\subsection{Concluding Remarks}

We investigated a system consisting of agents subject to exogenous constant inputs  and agents without any inputs. 
Future research will consider extensions to time-varying exogenous inputs and more general graph topologies. 

\vspace{0.25cm}


\bibliographystyle{IEEEtran}  
\baselineskip 12pt
\begin{thebibliography}{100}

\bibitem{ref:1}
M. Mesbahi and M. Egerstedt, ``Graph Theoretic Methods in Multiagent Networks,'' \textit{Princeton University Press}, 2010. 

\bibitem{ref:2}
R. Olfati-Saber, J. Fax, and R. M. Murray, ``Consensus and Cooperation in Networked Multiagent Systems,'' \textit{Proc. of the IEEE}, vol. 95, pp. 215--233, 2007. 

\end{thebibliography}
\end{document}